\documentclass[12pt]{amsart}

\usepackage{soul}

\def\cal{\relax}

\usepackage{DLdef1}
\let\ssec\subsection

\renewcommand {\ssbegin}[1]
 {\refstepcounter{subsection}
 \def \secno {\gdef \secno {}{\ssecfont
\thesubsection.\hskip 2ex}%
 }%
 \begin{#1}}

 \def\mmat #1,#2,#3,#4,{\text{\small\arraycolsep=3pt $
\begin{pmatrix}#1&#2\\#3&#4\end{pmatrix}$}}

\begin{document}

\markboth{Sofiane Bouarroudj, Pavel Grozman, Dimitry Leites, Irina
Shchepochkina} {Real-complex supermanifolds}

\thispagestyle{empty}

\title[Real-complex supermanifolds]{Minkowski superspaces and superstrings
as almost real-complex supermanifolds}

\author{SOFIANE BOUARROUDJ}

\address{New York University Abu Dhabi,
Division of Science and Mathematics, P.O. Box 129188, United Arab
Emirates; sofiane.bouarroudj@nyu.edu}

\author{PAVEL GROZMAN}

\address{Equa Simulation AB, Stockholm, Sweden; pavel.grozman@bredband.net}

\author{DIMITRY LEITES}

\address{Department of Mathematics, Stockholm University, Roslagsv. 101,
Kr\"aftriket hus 6, SE-106 91 Stockholm, Sweden; mleites@math.su.se}

\author{IRINA SHCHEPOCHKINA}

\address{Independent University of Moscow,
B. Vlasievsky per., d.  11, RU-119 002 Moscow, Russia;
irina@mccme.ru}




\begin{abstract} In 1996/7, J.~Bernstein observed that smooth or analytic supermanifolds
that mathematicians study are real or (almost) complex ones, while
Minkowski superspaces are completely different objects. They are
what we call \emph{almost real-complex supermanifolds}, i.e., real
supermanifolds with  a \emph{non-integrable} distribution, the
collection of subspaces of the tangent space, and in every subspace
a complext structure is given.

An almost complex structure on a real supermanifold can be given by
an even or odd operator; it is complex (without \lq\lq always") if
the suitable superization of the Nijenhuis tensor vanishes. On
almost real-complex supermanifolds, we define the circumcised analog
of the Nijenhuis tensor. We compute it for the Minkowski superspaces
and superstrings. The space of values of the circumcised Nijenhuis
tensor splits into (indecomposable, generally) components whose
irreducible constituents are similar to those of Riemann or Penrose
tensors. The Nijenhuis tensor vanishes identically only on
superstrings of superdimension $1|1$ and, besides, the superstring
is endowed with a contact structure. We also prove that all real
forms of complex Grassmann algebras are isomorphic although singled
out by manifestly different anti-involutions.
\end{abstract}


\thanks{Thanks are due to J.~Bernstein who made us interested in his
problem \lq\lq describe invariants of real-complex supermanifolds",
to A.~Lebedev for help and the referee for helpful comments. S.B.
was partly supported by the grant AD 065 NYUAD}



\maketitle

\section{Introduction}

\ssec{General remark on superizations}  In the supermanifold theory,
there are several \lq\lq straightforward" superizations of the
classical non-super notions. Definitions of superschemes and
supervarieties over any field  (\cite{Le0}), of their $C^\infty$
analogs --- supermanifolds (\cite{L1}), and of complex analytic
analogs --- superspaces (\cite{Va1}) are examples of such \lq\lq
straightforward" superizations. To be just, observe that to figure
out which direction of superization should be considered \lq\lq
straightforward" sometimes took a while, some of the above notions
are subjects of disputes and \lq\lq improvements" for more than 40
years by now.

There are also notions of Linear Algebra and algebraic notions of
Algebraic and Dif\-fe\-rential Geometries that have \emph{several}
superizations. Some of these superizations were (and some of them
still remain) unexpected and without direct non-super analog. For
example, among superdeterminants, there is the well-known Berezinian
$\Ber$ and several not so well known analogs (e.g., the queer
determinant $\qet$, and the \lq\lq classical limits" of $\Ber$ and
$\qet$, see \cite{DSB}, p. 476).

Studying supersymmetries may sometimes help not only to better
understand the classical non-super notion (like integral) but even
to distinguish a new notion in the non-super setting. Here we
consider one such notion, implicitly introduced together with
Minkowski superspaces. Lecturing on the results of this paper we
heard from the listeners that a CR-structure looks similar, and
indeed it does to an extent, but a careful comparison immediately
reveals that these notions have nothing in common, essentially.

\ssec{A new notion: Real-complex supermanifold} During the Special
year (1996--97) devoted by IAS, Princeton, to attempts to understand
at least some of the mathematics used in physical papers on
supersymmetry, J.~Bernstein pointed at one more example of an
unexpected super structure  (see notes of Bernstein's lectures taken
by Deligne and Morgan \cite{Del}, p. 94). It dawned upon him that
the models of our space-time (Minkowski superspaces) suggested in
the physical papers of pioneers, where supersymmetry was discovered,
are neither real nor complex supermanifolds, nor real supermanifolds
with a(n almost) complex structure; Minkowski superspaces are
different from real or complex supermanifolds or real supermanifolds
with a(n almost) complex structure studied by mathematicians so far
(e.g., see \cite{Va1,MaG}).

The Minkowski superspaces and superstrings introduced by physicists
are objects with a structure previously never considered. We give a
precise definition of such objects in the next subsection. Meanwhile
observe that although the bilinear forms (with Lorentzian signature)
given at the tangent space at every point of the Minkowski space are
equivalent, there are, nevertheless, several types of Minkowski
spaces. These spaces differ, for example, by the Riemannian
tensors\footnote{The fixed terms \lq\lq Riemann tensor", \lq\lq
Nijenhuis tensor" denote, strictly speaking, not tensors but tensor
\emph{fields}. In what follows we have to carefully distinguish
tensors from tensor fields.} constructed from the metrics.

Every $N$-extended Minkowski superspace and certain of the \lq\lq
super Riemann surfaces" considered in String Theories is a {\bf
real} ma\-ni\-fold (Minkowski space or a Riemann surface,
respectively) rigged with the sheaf of functions with values in
$\Lambda_\Cee(s)$, the {\bf complex} Grassmann superalgebra with $s$
generators. This construction differs from the case of \lq\lq the
sheaf of complex-valued functions on a real manifold" ($s=0$) in two
ways:

1) Considering the sheaf of $\Cee$-valued functions on the real
(super)manifold we do not get anything new as compared with
considering the sheaf of $\Ree$-valued functions on the same
(super)manifold. Indeed, there is only one real structure on the
target space $\Cee$, and hence every $\Cee$-valued function $f$ can
be {\bf canonically}  represented in the form of a sum of a pair of
$\Ree$-valued functions: $f(x)=u(x)+iv(x)$. If $s>0$, there is no
{\bf canonical}, i.e., unique in some way, real structure.

A question arises: how many isomorphism classes of real forms of the
Grassmann algebra $\Lambda_\Cee(s)$, i.e., the space of values of
superfunctions at a given point, are there?

Obviously, there are {\it several} analogs of the complex
conjugation on the Grassmann superalgebra of very distinct shape
(mathematicians favor some of them, physicists favor other ones, see
\cite{B, MaG, Del}). The answer to the above question was not given
anywhere, as far as we know, except \cite{L1}, where it was borrowed
from the first arXiv version of this paper. We will prove that the
answer is as follows: all the real forms of $\Lambda_\Cee(s)$ are
isomorphic (for $s<\infty$), but there is no \emph{canonical} real
form if $s>0$.

2) One more peculiarity of $N$-extended Minkowski superspace
$\cM_N$, whose definition we recall in subsec. 3.3, is the presence
of a non-holonomic (i.e., non-integrable) distribution $\cM_N$ is
rigged with.

If $\cM_N$ were just Minkowski space $M$ with functions on it taking
values in $\Lambda_\Cee(2N)$, this would have meant that on the
purely odd subspace of the tangent space at each point of $\cM_N$
there is given a complex structure. Having singled out the integral
submanifold $\cI$ of the distribution of codimension $4|0$ on
$\cM_N$ we could have offered a pair for characterization of
$\cM_N$: the Riemannian tensor on $M$ and, on $\cI$, the Nijenhuis
tensor describing the measure of non-flatness of the complex
structure on $\cI$. There are many distributions determined by
purely odd subspaces of the tangent spaces at the points of a given
supermanifold; some are integrable, some are not. The snag is: each
of the distribution that determine the Minkowski superspaces $\cM_N$
are non-integrable, so no integral subsupermanifold $\cI$ exists,
and nobody knew how to define the analog of the Nijenhuis tensor in
such a situation.

\ssec{Real-complex supermanifold as a supermanifold with a
$G$-structure and non-integrable distributions} There are known
various examples when a tensor field of a given type is given and,
in the tangent space at a given point, a \lq\lq flat" shape of the
tensor of an equivalence class is selected. The possibility of
reducing the tensor field under consideration to the  selected
\lq\lq flat" shape in an infinitesimal neighborhood of the given
point depends on the obstructions, cocycles representing cohomology
we will describe shortly.

{\bf Examples}: (1) on a given vector space $V$ over $\Cee$, there
is just one equivalence class of non-degenerate symmetric bilinear
forms $g$; (2) same is true for the non-degenerate anti-symmetric
bilinear forms $\omega$; (3) on a given vector space $V$ of
dimension $2n$ over $\Ree$, there is one equivalence class of the
automorphism $J$ such that $J^2=-\id$. (4) There can occur several
equivalence classes of tensors of a given type, e.g., non-degenerate
symmetric bilinear forms $g$ over $\Ree$ are distinguished by their
signature.

If $V$ is a vector space over $\Cee$ or $\Ree$ endowed with a tensor
$T$ whose automorphism group is $G$, and $M$ is a manifold  such
that $\dim M=\dim V$ and endowed with a tensor field whose value at
each point is equivalent to $T$, then $M$ is said to be endowed with
a $G$-\emph{structure}.

In what follows we recall and superize the notion of \emph{structure
functions}, i.e., functions on the principal $G$-bundle with values
in certain Lie algebra cohomology, see \cite{St}. Superization of
this notion is immediate and obvious. Obstructions to ``flatness''
of the three $G$-structures mentioned in examples above are the
Riemann tensor for metrics, $d\omega$ for almost symplectic
structures, and the Nijenhuis tensor for almost complex structures.

At every point of $N$-extended Minkowski superspace $\cM_N$, we can
select any shape of the tensor $J$ that defines the complex
structure on the odd subspace of the tangent space for a ``flat''
one. Question: what are the obstructions to reducing the tensor
field $J$ to the flat shape in an infinitesimal neighborhood of the
point? Here it is vital (nobody had ever considered such structures)
that \emph{an almost complex structure (tensor) $J$ is given on the
whole tangent spaces but only on the subspaces that constitute a
non-integrable distribution}.

The example of Minkowski superspaces can be naturally generalized.
Here is the most broad definition: an \emph{almost real-complex
supermanifold} of superdimension\footnote{In a recent preprint
\cite{Wi}, Witten used asterisk to separate real superdimension of
the underlying real supermanifold from that of the complex
superdimension of the superspaces of non-holonomic distribution in
particular cases $(p|0;*0|s)$. We find this notation very suggestive
and adopt it.} $(p|q;*r|s)$ is a real supermanifold
$\cE_\Ree^{p+2r|q+2s}$ endowed with a {\bf non-integrable}
distribution $\cD$ whose value at every point ---
$2r|2s$-dimensional subspace of the tangent space at the point ---
is endowed with a complex structure $J$.

In \cite{Del}, J.~Bernstein considered a particular case of
real-complex supermanifolds of superdimension $(p|0;*0|s)$ and used
a somewhat self-contradictory term \emph{cs manifolds}, short for
\emph{complex super manifolds} although these supermanifolds are not
complex.

Importance of non-integrability of the distribution $\cD$ is vital:
if $\cD$ were integrable, one could have restricted the problem onto
the integral subsupermanifold where the classical Nijenhuis tensor
provides with the answer.

In particular, almost real-complex structures exist on manifolds as
well, but only the example of Minkowski superspace drew attention to
them.

The tensor field $J$ that determines the structure of an almost
real-complex supermanifold, being defined on subspaces of the
non-holonomic distribution, defines a \emph{circumcised} (some say
\emph{reduced}) connection, see A.Vershik's appendix to \cite{Se}.
At the time J.Bernstein made his observation not only the definition
of the circumcised connection was corrected as compared with
\cite{Se}, the corresponding curvature tensor was known and even
computed for Minkowski superspaces, see \cite{GLs}, the details,
however, were published much later (see \cite{GL4}) and the
invariants of almost real-complex supermanifolds are described for
the first time.

\sssec{The new notion and superstrings} At the seminar on \lq\lq
super Riemann surfaces" and superstrings (for some of its workouts,
see \cite{LJ}), our attention was again drawn to a related problem
mentioned in \cite{Del}: Whereas all Riemann surfaces are
automatically endowed with a complex structure (and vice versa: each
complex curve possesses a naturally defined Riemannian metric), it
is completely unclear why should the analogs of these statements be
true for complex supercurves of superdimension $1|N$, and why real
supermanifolds of superdimension $2|2N$ must possess an almost
complex structure, to say nothing about complex one, even if the
underlying surface is endowed with a Riemannian metric.

As we will show, these great expectations, implicitly assumed in the
known to us texts on super analogs of Riemann surfaces, are
unjustified: the appropriate analogs of the Nijenhuis tensor are,
generally, non-zero, except for cuperstrings of superdimension
$1|1$, and not arbitrary but only those endowed with a contact
structure. Since in (super)conformal field theories an important
ingredient is the integral over the moduli (super)space of complex
structures --- (super)Teichm\"uller space $\cT$, it is worth to have
in mind that the Nijenhuis tensor $\cN_J$ should vanish for each of
the tensors $J$ parameterized by $\cT$.

\sssec{Question} What are the obstructions to integrability of the
almost complex structure on real supermanifolds and of the almost
real-complex structure on real-complex supermanifolds? In this paper
we answer this question and give several examples pertaining to
currently most popular physical models.

\ssec{Integrability of almost complex structures: Analytic and
algebraic approaches}

\sssec{Definitions} The real supermanifold of superdimension $2p|2q$
endowed with an even automorphism $J$ of the tangent space at every
point, i.e., a tensor field $J$ of valency $(1,1)$ such that
$J^2=-\id$ is said to be an {\it almost complex} supermanifold.

The supermanifold of superdimension $n|n$ over the ground field
$\Kee$ endowed with an odd tensor field $J$ of valency $(1,1)$ such
that $J^2=-\id$ is said to be an {\it almost $J$-symmetric}
supermanifold.

The supermanifold of superdimension $n|n$ over the ground field
$\Kee$ endowed with an odd tensor field $\Pi$ of valency $(1,1)$
such that $\Pi^2=\id$ is said to be an {\it almost $\Pi$-symmetric}
supermanifold.

Over $\Cee$, and also if the characteristic of the ground field
$\Kee$ is equal to 2, every almost $\Pi$-symmetric structure is
isomorphic to an almost $J$-symmetric structure, and the other way
round; accordingly,  the Lie supergroup $G=GQ_J(n)$ preserving the
odd operator $J$ such that $J^2=-\id$ is isomorphic to the
supergroup $G=GQ_\Pi(n)$ preserving the odd operator $\Pi$ such that
$\Pi^2=\id$. Contrariwise, the real supergroups $GQ_J(n; \Ree)$ and
$G=GQ_\Pi(n; \Ree)$ are not isomorphic, see \cite{L1}.

A given supermanifold $\cM$ is said to be \emph{endowed with an
almost $G$-structure} (the adjective \lq\lq almost" is often not
mentioned causing confusion) if, in the tangent space to every point
of $\cM$, an action of a supergroup $G$ is given (as in examples
above, where $G=GL(p|q;\Cee)\subset GL(2p|2q;\Ree)$ for $J$ even and
$G=GQ_J(n;\Kee)$ for $J$ odd). If this $G$-structure is integrable
(for details, see \cite{St, GL4}), one can drop the adjective \lq\lq
almost".

There are two ways to study obstructions to integrability, we expose
them in the next two subsections:

\sssec{Analytic approach} The obstruction to integrability of the
tensor field $J$ on ma\-ni\-folds were first, as far as we know,
computed by Newlander and Nirenberg \cite{NN}. The obstructions
constitute what is called the \emph{curvature} of the almost complex
structure given by $J$ or the \emph{Nijenhuis tensor} (field)
$\cN_J$ defined, for any vector fields $X, Y$, to be:
\begin{equation}\label{nijenhev}
\cN_J(X, Y)=[J(X), J(Y)]-J([J(X), Y]-J([X, J(Y)])-[X, Y],
\end{equation}
where $[\cdot, \cdot]$ is the Lie bracket of vector fields.
Nijenhuis showed that the tensor field $\cN_J$ can be represented in
the following form:
\begin{equation}\label{nijenhev1}
\cN_J=\nfrac12\{J, J\}_N,
\end{equation}
where $\{\cdot, \cdot\}_N$ is the {\em Nijenhuis bracket}, see
\cite{GrI}, where there are listed all bilinear differential
operators acting in the spaces of tensor fields and invariant with
respect to the changes of variables; see also \cite{Nij}, where 12
equivalent definitions of the Nijenhuis tensor are given.

A.~Vaintrob studied integrability of the \emph{even} tensor $J$ on
supermanifolds and  proved that an almost complex supermanifold is
complex if the straightforward super analog of the Nijenhuis tensor
given by the same expression \eqref{nijenhev} vanishes, see
\cite{Va2}. This fact (rediscovered in \cite{McH}) is used in a rich
with results paper \cite{Va1}.

The same applies as well to the $J$-symmetry, called in \cite{Va2}
\emph{odd complex structure}. The corresponding \emph{Nijenhuis
tensor} $\cN_J$ is defined, for any vector fields $X, Y$, to be:
\begin{equation}\label{nijenhod}
\cN_J(X, Y)=(-1)^{p(X)}[J(X), J(Y)]-J([J(X), Y])-(-1)^{p(X)}J([X,
J(Y)])-[X, Y].
\end{equation}
For any $\Pi$-symmetry, the corresponding tensor $\cN_\Pi$ is
defined by the same expression \eqref{nijenhod} with $J$ replaced by
$\Pi$.

So far, nobody bothered to investigate if there are other, apart
from $2|0$, superdimensions $p|q$ over $\Ree$ in which the Nijenhuis
tensor on $\cM^{p|q}$ vanishes identically, as in the case of almost
complex curves (=Riemannian surfaces). This negligence is
understandable because this task is difficult in the analytic
approach.

Here we show that the Nijenhuis tensor $\cN_J$ is irreducible
whereas the circumcised Nijenhuis tensors on Minkowski superspaces
and on various types of superstrings split into several components,
similar to $\alpha$- and $\beta$-components of the Penrose tensor,
cf. \cite{Po}. The Nijenhuis tensor on $\cM^{p|q}$ vanishes
identically if $p|q=1|1$ and $\cM^{1|1}$ is endowed with an almost
$J$-symmetry or $\Pi$-symmetry, and also if the $1|1$-dimensional
over $\Cee$ superstring $\cM^{1|1}$ is endowed with a contact
distribution and a circumcised almost complex structure on the folia
of the distribution (as in the Neveu-Schwarz and Ramond cases).

These discoveries were possible since there is, fortunately, an
approach more adequate for the task:

\sssec{Algebraic approach}\label{san} The same problems can be
formulated in terms of obstructions to integrability of
$G$-structures; compare the definitions in \cite{St} and \cite{GL4}.
In these terms superization is immediate and obvious. For example,
investigation flatness of the almost complex structure given by a
tensor $J$ we consider first the Lie superalgebra
$\fg_0:=\fgl(p|q;\Cee)$ preserving the tensor $J$ considered at a
given point, identify the tangent space at this point with the
$\fg_0$-module $\fg_{-1}$,  and construct the Cartan prolong
$\fg_*:=(\fg_{-1}; \fg_0)_*$ considered as a {\bf real} Lie
superalgebra (for the definition of Cartan prolong and its
generalizations we need below, see \cite{Shch} and Appendix). By the
usual arguments applicable to any $G$-structure (\cite{St}), the
obstruction to integrability of the almost $G$-structure is a
tensor, called \emph{structure function}, whose values at every
point lie in the space
\begin{equation}\label{curvtens}
H^2(\fg_{-1}; \fg_*).
\end{equation}
The curvature tensor itself is embodied by a cocycle representing a
non-trivial cohomology class of $H^2(\fg_{-1}; \fg_*)$.

For example, $H^2(\fg_{-1}; \fg_*)$ is the space of values (at the
point considered) of either the Riemann tensor if $\fg_0=\fo(n)$, or
of the Nijenhuis tensor if $\fg_0=\fgl(2n;\Ree)$, or of obstructions
to integrability of the almost symplectic structure if
$\fg_0=\fsp(2n)$, etc.

\paragraph{Example} Consider one example in more details.
On a manifold $M$, let there be given a non-degenerate
anti-symmetric bilinear form $B_m$  at each tangent space $T_mM$;
let $\omega$ be the exterior 2-form determined by a collection of
the forms $B_m$ for all points $m\in M$.

{\bf Definition}. Manifold $M$ with the above form $\omega$ is said
to be \emph{almost symplectic}. If the form $\omega$ can be reduced
to the chosen flat shape $\sum dp_i\wedge dq_i$ not only at every
point (which is always possible), but also in its (infinitesimal)
neighborhood, then the manifold is said to be \emph{symplectic}.

Let $V$ be a space isomorphic to $T_mM$, and $\fsp(V)$ the
symplectic Lie algebra preserving the image $B$ of the form $B_m$
under this isomorphism. As is known from Linear Algebra, the form
$B$ determines a canonical isomorphism $V\simeq V^*$, so to every
element $c\in \Hom(V\wedge V, V)$ we can assign an exterior 3-form:
\[C(u,v,w)=B(c(u,v),w)+ B(c(v,w), u)+B(c(w,u),v).\]
The above described map $c\mapsto C$ sends the coboundaries, i.e.,
elements $c$ of the form
$$
c(u,v)=S(u)v-S(v)u,\text{  where } S\in \Hom(V,\fsp(V)), \text{ and}
u,v\in V,
$$
to $0$, whereas the obstruction to flatness of $\omega$ at $m$ is
precisely the 3-form $C$ on $V$. The collection of all forms $C$ for
all points $m\in M$ constitute the exterior form $d\omega$, the
analog of the Rieamann tensor (field) for the almost symplectic
case. If $d\omega=0$, then $\omega$ can be reduced to the flat form
(Darboux's theorem).

\paragraph{How to express the obstructions to flatness}
It is convenient to represent the space $H^2(\fg_{-1}; \fg_*)$ as
the sum of irreducible $\fg_0$-modules (as this was done long ago
for the Riemann tensor represented as the sum of Weyl tensor,
traceless Ricci tensor and scalar curvature), and further shorthand
this representation by considering only highest weights of the
irreducible modules. Such shorthand expression is only possible if
the complexification of $\fg_0$ is either a semi-simple Lie algebra
or a central extension thereof, i.e. is a
{\it reductive complex Lie algebra}
If only $(\fg_0)_\ev$ is a reductive complex Lie algebra, then at
least one can describe how the irreducible $\fg_0$-modules are
glued; for examples of such descriptions, see \cite{Po} and further
in this paper.

\sssec{On analogs of Wess-Zumino constraints} The $\Zee$-grading of
the Lie algebra $\fg_*$ induces a $\Zee$-grading on the space
$H^2(\fg_{-}; \fg_*)$, called the \emph{degree}. As is shown in
\cite{St}, structure functions of a given degree are only defined
provided all the structure functions of smaller degrees vanish; same
is true in the non-holonomic situation. In supergravity theory these
conditions (vanishing of the structure functions of lesser degrees)
are known as \emph{Wess-Zumino constraints}, see \cite{GL4}.

\ssec{Superstrings and complex structures} Let us  begin with
supermanifolds of real superdimension $2|2m$, perhaps endowed with
an additional structure (such as preserved, infinitesimally, by the
centerless Neveu-Schwarz or Ramond superalgebra). In the String
Theory, these supermanifolds are usually referred to as \emph{super
Riemann surfaces}. We know that there are not two but four infinite
series and several exceptional simple Lie superalgebras analogous to
the simple Lie algebra of vector fields on the circle (for the
classification, see \cite{GLS}), called {\it stringy superalgebras}.
Exactly ten of them are {\it distinguished}: They are simple, and
have non-trivial central extensions, and hence are particularly
interesting from the point of view of possible physical
applications, since only these centers might act on super versions
of certain infinite dimensional Teichm\"uller spaces (parameterizing
deformations of the complex structures), cf. \cite{Kon, BSh}. So we
should investigate integrability of the appropriate Nijenhuis tensor
on the supersurfaces rigged with one of the structures preserved by
the distinguished stringy superalgebra.

As on Minkowski superspaces, some of these structures (e.g., in the
centerless Neveu-Schwarz and Ramond cases) are {\it nonholonomic}
(i.e., non-integrable) distributions, like the contact one.

\ssec{Analogs of the curvature tensor in presence of non-integrable
distribution} Recall the definition of the analogs of the curvature
tensor in presence of non-integrable distribution (for more details,
see \cite{GL4}). First of all, in order to distinguish the
non-integrability of of the distribution from the
(non-)integrability of the curvature field, we will,  speaking about
distributions, say \emph{non-holonomic}. The arguments of \cite{GL4}
are literally superized, so in this subsection we drop \lq\lq
super".

Take the $\Zee$-graded Lie algebra $\fg_*:=(\fg_{-},
\fg_0)_*=\mathop{\oplus}\limits_{i\geq -d}\fg_i$ associated with the
filtered Lie algebra preserving the distribution and the tensor
$\cT$ whose flatness we are studying. Therefore, we identify the
tangent space at every point with the space
$\fg_{-}=\mathop{\oplus}\limits_{i<0}\fg_i$. If the Lie algebra
$\fg_{-}$ is generated by the space $\fg_{-1}$, then this space is
precisely the value of the distribution at the point, whereas
$\fg_{0}$ is the Lie algebra that preserves the tensor $\cT$.

By analogy with the representation of the Nijenhuis tensor (as well
as obstructions to flatness of any other $G$-structure) as a cocycle
representing a class from $H^2(\fg_{-1}; \fg_*)$, see
\eqref{curvtens}, we give the following

{\bf Definition}. The \emph{circumsized curvature tensor} of the
structure given by the tensor $\cT$ in presence of a non-holonomic
distribution $\cD$ is a cocycle representing a non-trivial class of
\begin{equation}\label{Ncurvtens}
H^2(\fg_{-}; \fg_*),\text{~~where
$\fg_{-}=\mathop{\oplus}\limits_{i<0}\fg_i$}.
\end{equation}

Let $\widetilde\fg_*:=(\fg_{-}, \widetilde\fg_0)_*$ be a simple Lie
algebra such that $\fg_0$ has center $\fz$ and the extension of
$\widetilde\fg_0=\fg_0/\fz$ to $\fg_0$ is trivial, i.e., the center
is a direct summand. In \cite{GL4}, it is shown that
$H^2(\fg_{-};\widetilde\fg_*)$, the space of values of
\emph{non-holonomic analog of the Riemannian tensor} strictly
contains $H^2(\fg_{-};\fg_*)$, the space of values of
\emph{non-holonomic analog of the Weyl tensor}. The case where the
Lie algebra of the Lie group $G$ that defines the $G$-structure is
$\fg_0$ (resp. $\widetilde\fg_0$) is said to be \emph{conformal}
(resp. \emph{reduced}).

\ssec{Our results} 1) We classified the real structures on the
finite dimensional Grassmann algebra $\Lambda_\Cee(s)$, see \S2. We
prove that, although the anti-automorphisms that single out the
\emph{real forms} of $\Lambda_\Cee(s)$ look completely differently,
the real forms they determine are isomorphic, albeit
non-canonically.

2) We described the two types of obstructions to integrability:

(2a) of the almost complex structure, as well as almost
$J$-symmetric and $\Pi$-symmetric structures, given on the whole
tangent space at a point (super versions of the Nijenhuis tensor),
and

(2b) of the almost real-complex structure given on the subspaces of
a non-holonomic distribution (a circumcised version of the Nijenhuis
tensor), see \S3.

Obstructions to integrability of an almost complex structure
constitute an irreducible $\fg_0$-module, while the obstructions to
integrability of an almost $J$-structure and $\Pi$-structure form an
indecomposable $\fg_0$-module, and we described how it is glued of
irreducibles.

Obstructions to integrability of an almost real-complex structure
are illustrated with examples of Minkowski superspaces and
superstrings. In particular, we show that the circumcised Nijenhuis
tensor identically vanishes on $1|1$-dimensional complex superstring
with contact structure, as predicted by P.~Deligne in 1987, see
\cite{MaT}. We show further that this is the only superdimension and
structure on the superstring (among the four infinite series and
four exceptional stringy superalgebras) for which the circumcised
Nijenhuis tensor vanishes identically. The Nijenhuis tensors of
$\Pi$- and $J$-symmetric (``odd complex'') structures vanish
identically only in dimension $1|1$.

Lorentzian metric on the Minkowski space $M$ underlying the
superspace $\cM_N$ is induced by the non-holonomic distribution of
totally even codimension, so the circumcised Nijenhuis tensor
replaces, for $N>0$, the Riemannian tensor on $M$.

The Minkowski space $M$ can be considered as the quotient of the
Poincar\'e group modulo Lorentz group, but can also be considered
\lq\lq conformally", as the twistor space, see \cite{MaG}. We
considered the Minkowski superspaces $\cM_N$ from super versions of
both these approaches.

The statements of \S3 are obtained by means of the \texttt{SuperLie}
package \cite{Gr} and can be verified by the usual cohomology
technique described in detail in \cite{Po}.

Minkowski space can possess lots of various metrics that differ by
their Riemann tensor --- the measure of their \lq\lq non-flatness".
Similarly, there are lots of different Minkowski superspaces and
superstrings that differ by circumcised and usual Nijenhuis tensors,
described in theorems \ref{Nij}, \ref{qpi}, \ref{mink},
\ref{string}.

3) In \S4 we showed that there are four ways to superize the notion
of K\"alerian manifold; there are even more ways to superize the
notion of a hyper-K\"alerian manifold. We also explained that these
superized notions can be endowed with almost real-complex
structures. All these superizations seem to be new.

\ssec{Related problems} Our methods are applicable to non-standard
models of complexified and compactified Minkowski superspaces, such
as the ones considered in \cite{MaG,GL}.

The Nijenhuis tensor is \lq\lq inexhaustible as atom", and recently
O.~Bogoyavlenskij described, on manifolds, several of its
properties previously unnoticed, see \cite{B}. It is interesting to
consider the super version of Bogoyavlenkij's problem.

{\bf Mixture of even and odd complex structures}. Observe that the
{\em odd complex structure} such as $J$-symmetry (as well as
$\Pi$-symmetry) can be defined over $\Cee$, so an almost complex
supermanifold might be also endowed with an odd almost complex
structure ($J$-symmetry or $\Pi$-symmetry). Such \lq\lq doubly
complex" structures exist on the superstrings with Ramond
superalgebra (or rather its quotient modulo center) as its Lie
superalgebra of infinitesimal symmetries, cf. \cite{LJ}, and these
even and odd complex structures are integrable or not independently.

Our approach to computing the analogs of the Nijenhuis tensor is
applicable to any of these open problems.

\section{All real forms of the Grassmann algebra are isomorphic}

\ssec{Real forms of the complex superalgebra} Given a superalgebra
$C$ over $\Cee$, we say that an even $\Ree$-linear map $\rho: C\tto
C$ is a {\it real structure} on~$C$ if
\begin{equation}
  \label{ab_rho_1}
\rho^{2} = \id,\quad   \rho(ab) = \rho(a) \rho(b) ,\text{~ and ~}
  \rho(za) = \bar{z}\rho(a) \quad
  \text{ for any $z \in \Cee$ and $a, b \in C$}.
\end{equation}

We set
\begin{equation}
  \label{Re_C__Im_C}
\begin{split} \RE_\rho C = \{ a~\in C \mid \rho(a)=a\} ,\quad
  \IM_\rho  C = \{ a~\in C \mid \rho(a) =-a \}.
  \end{split}
\end{equation}
Recall that the \emph{realification} $C^\Ree$ of $C$ is the same $C$
but considered over $\Ree$. Clearly, $\RE_\rho C$ is a subalgebra in
the realification $C^\Ree$ of $C$, and $\IM_\rho C = i\cdot \RE_\rho
C $ whereas $C^\Ree= \RE_\rho C \oplus \IM_\rho C$. The subalgebra
$\RE_\rho C$ is said to be a \emph{real form} of $C$ (corresponding
to the real structure $\rho$).

Observe that, on $C^\Ree$, the map $\rho$ is an automorphism, and
hence, if $C$ is an algebra with unit, then $\rho(1)=1$, so, $1\in
\RE_\rho C$.

\sssec{Examples of real structures on the Grassmann superalgebra
$\Lambda_\Cee(n)$} Let $\theta=(\theta_1, \dots \theta_n)$ be
generators of  $\Lambda_\Cee(n)$. For $m=0$, there is, obviously,
just one real structure, the canonical one.

For $n=1$, there are many real structures; clearly,
$\rho(\theta)=\lambda\theta$ for $\lambda\in\Cee$. Since $\rho$ is
involutive, $\lambda\bar\lambda=1$, i.e., $\lambda=\exp(i\varphi)$,
where $\varphi\in\Ree$; and hence
\[ {\RE_\rho
\Lambda_\Cee(1)=\{a+b\exp(i\varphi/2)\theta\mid a,b\in\Ree\}}.
\]

For $n=2k$, set $\theta=(\xi, \eta)$, where $\xi=(\xi_1, \dots
\xi_k)$, $\eta=(\eta_1, \dots \eta_k)$. The following are the main
examples of real structures; the first one favored by
mathematicians, the second one by physicists; one can (and for $n$
odd, one should) consider a mixture of these structures:

\begin{enumerate}
  \item $\rho_{bar}(\theta_j)=\theta_j$ for any $j=1, \dots , n$
  (obviously, one may consider an $n$-parameter generalization
  $\rho_{bar}(\theta_j)=\exp(i\varphi_j)\theta_j$);
  \item $\rho_{tr}(\xi_j)=i\cdot \eta_j$; $\rho_{tr}(\eta_j)=i\cdot
  \xi_j$ for any $j=1, \dots , k$ and $i=\sqrt{-1}$.
\end{enumerate}

In view of the above described diversity of involutions that single
out real forms of the complex Grassmann algebra (this diversity is
discussed at length in \cite{Ber} and mentioned in \cite{MaG} and
\cite{Del}), the following theorem, although obvious to some
experts, is worth mentioning.  To the best of our knowledge it was
not even formulated so far.

\ssbegin{Theorem}\label{relFofGr} All real forms of the Grassmann
algebra are isomorphic.\end{Theorem}

\begin{proof} Let $G:=\Lambda_\Cee(n)=\Lambda_\Cee(\theta)$, where
$\theta=(\theta_1, \dots , \theta_n)$ are generators.  Set
\[
G_k=\mathop{\oplus}\limits_{s\ge k} \Lambda^{s}(\theta),\
G_\ev=\mathop{\oplus}\limits_{2s\le n} \Lambda^{2s}(\theta),\
G_\od=\mathop{\oplus}\limits_{2s-1\le n} \Lambda^{2s-1}(\theta).
\]
For $n$ odd, introduce also the space
$G_\od^-=\mathop{\oplus}\limits_{1\le 2k-1<n} \Lambda^{2k-1}$.

Let $\rho$ be a real structure on $G=\Lambda_\Cee(n)$. Observe that
if $U\subset G$ is a  $\rho$-invariant complex linear space, i.e.,
$\rho(U)=U$, then $U$, same as $G$, splits into the direct sum of
its real and imaginary parts: $U^\Ree= \RE_\rho U \oplus \IM_\rho
U$, and $(G/U)^\Ree= \RE_\rho G/U \oplus \IM_\rho G/U$  if $U$ is an
ideal.

Having selected anticommuting generators $\theta=(\theta_1, \dots ,
\theta_n)$ of $G$, we construct the ideals $G_k$. For the natural
filtration of the Grassmann superalgebra $G=G_0\supset G_1\supset
\dots\supset G_n$ associated with the degree of the elements of $G$
assuming that each generator is of degree 1, set $V=G_1/G_2$, and
let $\pi:G_1\tto V$ be the natural projection.

Observe that

1) $G=\Cee\cdot 1\oplus G_1$, as linear space;

2) $G_\od\subset G_1$.

Now observe that $G_1$ has an invariant description: this is the
ideal of nilpotent elements of $G$. Therefore $\rho(G_1)=G_1$. But
then $\rho(G_k)=G_k$ for all $k=1,\dots, n$.

Observe that the center of $G$ is
\[Z=\begin{cases}
G_\ev&\text{for $n$ even}\\
G_\ev\oplus G_n&\text{for $n$ odd}.\end{cases}\] For any $n$, we
have $\rho(Z)=Z$ and $\rho(Z\cap G_2)=Z\cap G_2$.
For $n$ odd, set
$G_\od^-=\mathop{\oplus}\limits_{1\le 2k-1<n} \Lambda^{2k-1}$.

Now, let $B$ be the real form of $G$ corresponding to the real
structure $\rho$, let $B_k$ be the real form of the ideal $G_k$, and
$B_Z$ the real form of $Z\cap G_2$. Then $\pi(B_1)$ is the real form
of the space $V$, and hence $\dim_\Ree \pi(B_1)=n$.

Let $x_1, \dots, x_n\in B_1$ be such that $\pi(x_1),\dots \pi(x_n)$
form a basis of $V$. Clearly, $x_1,\dots, x_n$ generate the algebra
$B$ over $\Ree$ and the algebra $G$ over $\Cee$. Let us expand each
$x_k$ with respect to $Z\cap G_2$ and $G_\od$ for $n$ even and with
respect to $Z\cap G_2$, $G_\od^-$ and $G_n$ for $n$ odd:
\[
x_k=y_k+z_k, \text{ ~where~ } y_k\in \begin{cases}G_\od&\text{for
$n$ even},\\G_\od^-&\text{for $n$ odd},\end{cases} \; z_k\in Z\cap
G_2.
\]

Since $z_k \in G_2$, it follows that $\pi(y_k)=\pi(x_k)$. Therefore
the elements $y_1,\dots, y_n$ also generate $G$ (over $\Cee$), and
anti-commute since they belong to $G_\od$. Hence so do their images
$\rho(y_1),\dots, \rho(y_n)$.

Now observe that since $x_k\in B$, it follows that $\rho(x_k)=x_k$,
i.e., \[\rho(y_k)+\rho(z_k)= y_k+z_k.\]  But $\rho(z_k)\in Z\cap
G_2$ implying that $\rho(y_k)=y_k+z'_k$, where $z'_k\in Z\cap G_2$.
Since the $\rho(y_1),\dots, \rho(y_n)$ anticommute, it follows that,
for all $k,l$, we have:
\begin{equation}\label{*}
y_kz'_l+z'_ky_l+z'_kz'_l=0 \Longrightarrow y_kz'_l+z'_ky_l=0 \text{
and }z'_kz'_l=0.
\end{equation}
The second and the third equalities follow from the fact that the
first two summands in the first equality lie in $G_\od$, and the
third summand lies in $G_\ev$.

Consider now the elements \[t_k=\frac 12(y_k+\rho(y_k))=y_k+ \frac
12z'_k\in B.\] Since $\pi(t_k)=\pi(y_k)$, the elements $t_k$
generate the algebra $G$ over $\Cee$ and the algebra $B$ over
$\Ree$. Thanks to \eqref{*} the elements $t_k$ anticommute.

The theorem now follows from the universality of the Grassmann
algebra as algebra with unit and anti-commuting generators.
Universality is understood here in the sense that any other algebra
with anti-commuting generators is a quotient of the Grassmann
algebra (just because there are no other relations in the Grassmann
algebra). Since we have found $n$ anti-commuting generators of the
algebra $B$ over $\Ree$, it follows from the dimension
considerations that there are no relations which are not corollaries
of anti-commutation ones.
\end{proof}

\section{Obstructions to integrability of almost complex and almost real-complex
structures}

\ssec{Revision of the classical examples} {\bf 1) The
$2n|2m$-dimensional almost complex supermanifolds with $J$ even}. It
is well-known that any almost complex structure on the real
orientable surface is integrable, see \cite{NN}. Let us investigate
if there are other exceptional cases where a given almost complex
structure is always complex, except superdimension $2|0$, and
investigate if the space of values of the Nijenhuis tensor can be
split into the irreducible $\fgl(n|m;\Cee)^\Ree$-modules.
(E.Poletaeva performed similar calculations for the analogs of the
Riemann tensor: There is no complete reducibility, and the
description of how the irreducible components are glued together is
rather intricate, see \cite{Po}.)

In terms of subsec. \ref{san} we consider the Lie superalgebra
$\fg_0=\fgl(n|m;\Cee)^\Ree$ consisting of supermatrices in the
non-standard format of the form
\begin{equation}\label{ab-ba}
\begin{pmatrix}
A&B\\
-B&A
\end{pmatrix},\text{~~where $A, B\in\fgl(n|m;\Ree)$, and $A+iB\in\fgl(n|m;\Cee)$
 for
$i=\sqrt{-1}$,}
\end{equation} and the tautological $\fg_0$-module
$\fg_{-1}=\Ree^{2n|2m}$ with the following format of its basis
vectors (even $|$ odd  $|$ even $|$ odd):
\[\Span(\partial_1,\dots, \partial_{n}\mid
\partial_{n+1},\dots, \partial_{n+m}\mid\partial_{n+m+1},\dots, \partial_{2n+m}\mid
\partial_{2n+m+1},\dots, \partial_{2n+2m}).\]

The classical results on manifolds state that the tensor field
$\cN_J$, see \eqref{nijenhev}, is the only obstruction to
integrability of the almost complex structure on manifolds; besides,
on surfaces, the tensor field $\cN_J$ vanishes identically.
According to \cite{Va2, McH}, on supermanifolds, $\cN_J$, see
\eqref{nijenhev}, is also the only obstruction to integrability of
the almost complex structure.

According to \cite{Va2}, the only obstruction to integrability of
the almost $J$-symmetry is $\cN_J$, see \eqref{nijenhod}.

Let us sharpen these claims. For $n+m>1$, we consider the natural
division of root vectors with respect to which the positive ones are
those above the diagonals of $A$ and $B$ in \eqref{ab-ba}.

\sssbegin{Theorem}\label{Nij} The lowest weight (with respect to the
$\fg_0$-action) cocycles representing the elements of $H^2(\fg_{-1};
\fg_*)$ for $\fg_0=\fgl(n|m;\Cee)^\Ree$ and its tautological module
$\fg_{-1}$ are as follows  (all of degree $1$):
\begin{equation}\label{n,0}
\begin{array}{lll}
n=1|m=0:& 0&\\
n>1|m=0:&\partial_1  \otimes  (\partial_{n-1}^* \wedge
\partial_{2 n}^*-\partial_n^* \wedge \partial_{2n-1}^*
), & \partial_{n+1} \otimes (\partial_{n-1}^* \wedge \partial_{2
n}^*
-\partial_n^* \wedge \partial_{2n-1}^* )\\
n=0|m\geq 1:&\partial_1\otimes (\partial_{m}^*\wedge \partial_{2m}^*
),&
\partial_1\otimes (\partial_m^*)^{\wedge 2}\\
n>0|m>0:&\partial_1\otimes (\partial_{n+m}^* \wedge \partial_{2m+2n}^*),& 
\partial_1 \otimes (\partial_{n+m}^*)^{\wedge 2}\\
\end{array}
\end{equation}

\end{Theorem}

\parbegin{Comment}\label{comment} 1) For $n|m=0|1$, the cocycles
\eqref{n,0} span the whole space $H^2(\fg_{-1}; \fg_*)$, not just
that of lowest weight vectors.

2) Since the space $H^2(\fg_{-1}; \fg_*)$ has two $\fg_0$-lowest
weight vectors, see \eqref{n,0}, the researcher familiar with
representations of Lie algebras over $\Cee$ might think that the
Nijenhuis tensor splits into two irreducible components. This is not
so since the ground field is $\Ree$: Each lowest weight vector is
obtained from the other one by multiplication by $i$ represented by
the matrix $\begin{pmatrix}
0&1_{n+m}\\
-1_{n+m}&0
\end{pmatrix}$,
and therefore the $\fg_0$-module $H^2(\fg_{-1}; \fg_*)$ is
irreducible; it is the realification of an irreducible module over
$\Cee$.
\end{Comment}

\sssec{Almost $J$-symmetric and almost $\Pi$-symmetric
supermanifolds over any ground field $\Kee$} Let $A,
B\in\fgl(n;\Kee)$, and $p(A)=\ev$, $p(B)=\od$; let
$\fg_0=\fq_\Pi(n;\Kee)$ (or $\fg_0=\fq_J(n;\Kee)$)be the Lie
superalgebra consisting of supermatrices of the form
\begin{equation}\label{abba}
\begin{pmatrix}
A&B\\
B&A
\end{pmatrix}\text{~~if $\fg_0=\fq_\Pi(n;\Kee)$, or }
\begin{pmatrix}
A&B\\
-B&A
\end{pmatrix}\text{~~if $\fg_0=\fq_J(n;\Kee)$},
\end{equation}
and $\fg_{-1}=\Kee^{n|n}$ be the tautological $\fg_0$-module with
the following basis vectors in the standard format (even $|$ odd):
\[\Span(\partial_1,\dots, \partial_{n}\mid
\partial_{n+1},\dots, \partial_{2n}).\]

\sssbegin{Theorem}\label{qpi} The cocycles representing the elements
of $H^2(\fg_{-1}; \fg_*)$ are all of degree $2$. The $\fg_0$-module
$H^2(\fg_{-1}; \fg_*)$ is indecomposable with the following
$(\fg_0)_\ev$-highest weights:
\begin{equation}\label{q}
\begin{array}{ll}
n=1:& \text{{\bf the zero module};
}\\
n=2:&\text{$(2,0)$, and $(1,1)$ each of multiplicity $2|2$ ($2$ even
ones and $2$ odd ones);
}\\
&\text{contains a submodule with the $(\fg_0)_\ev$-highest weights
$(2,0)$, and $(1,1)$}\\
&\text{each of multiplicity $1|1$;
}\\

n\geq 3:&\text{$\varphi_1:=(2,1,0,\dots,
0, -1)$, $\varphi_2:=(2,0,\dots, 0)$, and $\varphi_3:=(1,1,0, \dots, 0)$}\\
&\text{each of multiplicity $2|2$;
}\\
&\text{contains a submodule with the $(\fg_0)_\ev$-highest weights
$\varphi_1$ of multiplicity $2|2$, }\\
&\text{$\varphi_2$, and $\varphi_3$ each of multiplicity $1|1$;
}\\
&\text{containing, in turn, a submodule with the
$(\fg_0)_\ev$-highest weights $\varphi_2$, and $\varphi_3$}\\
&\text{each of multiplicity $1|1$.
}\\
\end{array}
\end{equation}

\end{Theorem}

\parbegin{Comment}\label{commentQ} 1) We were able
to shorthand the answer in Theorem \ref{qpi} by using the fact that
the superdimension of the space of vacuum (in our case, lowest
weight) vectors of every irreducible finite dimensional
$\fq_\Pi(n;\Kee)$- and $\fq_J(n;\Kee)$-module is of the form $k|k$
(except for the trivial module $\Kee$ and $\Pi(\Kee)$, when it is
equal to $1|0$ and $0|1$, respectively). Therefore, it suffices to
describe only half of the highest weight vectors, say the even ones.

2) The explicit answer is, however, hardly needed in theoretical
constructions neither in this theorem nor in theorem \ref{mink};
important is that the Nijenhuis tensor does or does not vanish
identically, and if there are several components what is the meaning
of vanishing of some of them: compare with the Einstein equations
$=$ vanishing of certain components of the Riemann tensor.
\end{Comment}

\ssec{The $N$-extended Minkowski superspaces $\cM_N$} Usually,
mathematicians represent complexified Minkowski superspaces as
homogeneous superspaces $\cM_N^\Cee=G_N^\Cee/P_N^\Cee$, where
$G_N^\Cee$ and $P_N^\Cee$ are certain complex supergroups. In the
papers and books written by physicists only the structure of the
tangent space $T_m\cM_N$ at a given point $m$ of Minkowski
superspace $\cM_N$, and the Lie superalgebra acting on $T_m\cM_N$
are usually given, actually, cf. \cite{GIOS} and references therein.

If  $\fg={\rm Lie}(G_N)$ and  $\fp={\rm Lie}(P_N)$, we have the
following: $\fg$ consists of supermatrices of the form
\begin{equation}\label{M1}
\fg=\Span\left(\begin{pmatrix}
A&0&0\\
Q&B&0\\
T&-\overline Q^t&-\overline A^t
\end{pmatrix},\;\;\begin{array}{l}\text{where $A\in\fsl(2;\Cee)$,
$Q\in\Mat_\Cee(2\times N)$,}\\ \text{$T=\overline T^t$,
$B\in\fgl(N;\Cee)$}\end{array}\right)
\end{equation} where bar denotes component-wise
complex conjugation, the $\Zee$-grading is block diagonal-wise
($\deg A=\deg B=0$, $\deg Q=-1$, $\deg T=-2$) and the parity is
defined as $\deg\mod 2$.

Let the Lie superalgebra $\fp$ consist of supermatrices of degree
$\ge 0$. Then $T_m\cM_N$ is spanned over $\Ree$ by matrices $T$
(constituting $T_mM$) and odd matrices $Q$.

The presence of a non-holonomic distribution on $\cM_N$ is
obvious\footnote{Recall a criterion of integrability of
distributions
--- {\bf  Frobenius's theorem}: \emph{A distribution is integrable if and only
if the sections of the distribution form a Lie sub(super)algebra
relative the bracket of vector fields}.}: the subspaces of the
tangent spaces that constitute the distribution are spanned by
vectors realized by pair of matrices $Q$ and $-\overline Q^t$, while
the superbracket of such vectors does belong to the linear
combination not of them but of the matrices $T$.

We can as well consider the maximal symmetry supergroup $G_N$ of
$\cM_N$, i.e., assume that $\fg=(\fg_-,\fg_0)_*$, where
$\fg_0=\fder_0(\fg_-)$ consists of the grading-preserving
derivations of the Lie superalgebra $\fg_-$. Set
$\fp:=\mathop{\oplus}\limits_{i\geq 0}\fg_i$. This approach and the
natural (especially in the light of successes of twistor models)
desire to have a {\bf simple} Lie superalgebra as the generalized
prolong, or its complexification, imposes restrictions on $B$ and
$A$, tying the group of inner symmetries with the Lorentz group.

Namely, this desire is satisfied if in the realization of the Lie
superalgebra by supermatrices of the form \eqref{M1} we ser
\begin{equation}\label{Mtw}
B\in\fu(N),\text{~ and $\tr B=\tr(A-\overline A^t)$, where
$A\in\fgl(2;\Cee)$.}
\end{equation}
Thus, the desire to have as the generalized Cartan prolong
a\emph{simple} Lie superalgebra with the same negative part as
\eqref{M1}, \emph{forces} us to consider on the Minkowski space $M$
underlying $\cM_N$ not a Lorentzian metric preserved by elements
$A\in\fsl(2;\Cee)\simeq \fo(3,1)$, but rather a \emph{conformal}
structure that only preserves the conformal class of the metric,
allowing to multiply the metric by non-zero constants. Then
$\fg_*=(\fg_-,\fg_0)_*$ consists of supermatrices of the form
\begin{equation}\label{M3}
\begin{pmatrix}
A&-\overline R^t&U\\
Q&B&R\\
T&-\overline Q^t&-\overline A^t
\end{pmatrix}\begin{array}{l}\text{~~with $Q$, $T$ as in
\eqref{M1}, $B=-\overline B^t$,
$R\in\Mat_\Cee(2\times N)$,}\\
\text{~~$\tr B=\tr(A-\overline A^t)$, $A\in\fgl(2;\Cee)$,
$U=\overline U^t$.}\end{array}
\end{equation}

\sssbegin{Theorem}\label{mink} In this theorem, the supermatrix
format is $2|N|2$ and $N=1$. Let $X_{i,j}$ stand  for the $(i,j)$-th
matrix unit in the matrix $X=A$, $B$, $T$, or $Q$.

{\bf 1) The \lq\lq conformal" case}. The Lie superalgebra
$\fg_*:=(\fg_-,\fg_0)_*$ consists of supermatrices of the form
\eqref{M3}.

Let the superscript denote the degree of the cocycle $c$, the
subscript its number. The highest weight cocycles representing the
basis elements of $H^2(\fg_{-1}; \fg_*)$ are as follows.
\[
 c_1^0 = T_{2,2}\otimes (Q_{1,1}^*
){}^{\wedge 2},\quad c_2^0 =T_{2,2}\otimes Q_{1,1}^* \wedge i
Q_{1,1}^*
\]

In degrees $1$ and $2$: None.

In degree $3$: (the numbering of cocycles match that of the reduced
case)
\[\tiny
\begin{array}{ll}
c_8^3 &=4Q_{1,2}\otimes (T_{1,1}^*\wedge
(T_{1,2}+T_{2,1})^*)+i(A_{1,1}-A_{2,2})\otimes (iQ_{1,1}^*\wedge
T_{1,1}^*)+i(A_{1,1}+A_{2,2})\otimes (iQ_{1,1}^*\wedge T_{1,1}^*)-\\
&4A_{1,2}\otimes (Q_{1,1}^*\wedge (T_{1,2}+T_{2,1})^*)+2
A_{1,2}\otimes (iQ_{1,1}^*\wedge i(T_{1,2}-T_{2,1})^*)+ 4
A_{1,2}\otimes (Q_{1,2}^*\wedge T_{1,1}^*)+\\
&2iA_{1,2}\otimes (iQ_{1,1}^*\wedge (T_{1,2}+T_{2,1})^*)- 2
iR_{1,1}\otimes (Q_{1,1}^*\wedge iQ_{1,1}^*)
+2 R_{1,1}\otimes (iQ_{1,1}^*)^{\wedge 2} \\[3mm]

 c_9^3 & = 4 iQ_{1,2}\otimes (T_{1,1}^*\wedge (T_{1,2}+T_{2,1})^*)-i(A_{1,1}-A_{2,2})\otimes (Q_{1,1}^*\wedge
 T_{1,1}^*)-i(A_{1,1}+A_{2,2})\otimes (Q_{1,1}^*\wedge T_{1,1}^*)-\\
 & 2 A_{1,2}\otimes (Q_{1,1}^*\wedge
 i(T_{1,2}-T_{2,1})^*)-4 A_{1,2}\otimes (iQ_{1,1}^*\wedge (T_{1,2}+T_{2,1})^*)+ 4 A_{1,2}\otimes (iQ_{1,2}^*\wedge
 T_{1,1}^*)-\\
 &2 iA_{1,2}\otimes (Q_{1,1}^*\wedge (T_{1,2}+T_{2,1})^*)+ 2 iR_{1,1}\otimes
 (Q_{1,1}^*)^{\wedge 2}-2 R_{1,1}\otimes (Q_{1,1}^*\wedge iQ_{1,1}^*)
\end{array}
\]

{\bf 2) The \lq\lq reduced" case}. Everything is as above but $B=0$,
$A\in\fsl(2;\Cee)$. In this case, $\fg_1=0$. The highest weight
cocycles representing the basis elements of $H^2(\fg_{-1}; \fg_*)$
are as in the conformal case with the following modifications or
additions: Two new cocycles in degree $1$ appear:
\[\footnotesize
\begin{array}{lcl}
c_3^1 & = & -Q_{1,1}\otimes (Q_{1,1}^*\wedge
iQ_{1,1}^*)+iQ_{1,1}\otimes (Q_{1,1}^*)^{\wedge 2}-Q_{1,2}\otimes
(Q_{1,1}^*\wedge iQ_{1,2}^*)+
iQ_{1,2}\otimes (Q_{1,1}^*\wedge Q_{1,2}^*)  \\[2mm]
   c_4^1& =&-Q_{1,1}\otimes (iQ_{1,1}^*)^{\wedge 2} +iQ_{1,1}\otimes
   (Q_{1,1}^*\wedge iQ_{1,1}^*)-Q_{1,2}\otimes (iQ_{1,1}^*\wedge iQ_{1,2}^*)
   +iQ_{1,2}\otimes (iQ_{1,1}^*\wedge Q_{1,2}^*)  \\[3mm]
\end{array}
\]

Three new cocycles in degree $2$ appear:

\[\tiny
\begin{array}{lcl}
c_5^2 & = & -i(A_{1,1}-A_{2,2})\otimes (Q_{1,1}^*)^{\wedge
2}-i(A_{1,1}-A_{2,2})\otimes (iQ_{1,1}^*)^{\wedge 2}+2
A_{1,2}\otimes (Q_{1,1}^*\wedge iQ_{1,2}^*)-\\&& 2 A_{1,2}\otimes
(iQ_{1,1}^*\wedge Q_{1,2}^*)-2 i A_{1,2}\otimes (Q_{1,1}^*\wedge
Q_{1,2}^*)-2 i A_{1,2}\otimes (iQ_{1,1}^*\wedge
iQ_{1,2}^*)-\\&&Q_{1,1}\otimes (iQ_{1,1}^*\wedge T_{1,1}^*)+
iQ_{1,1}\otimes (Q_{1,1}^*\wedge T_{1,1}^*)+ 2 Q_{1,2}\otimes
(Q_{1,1}^*\wedge i(T_{1,2}-T_{2,1})^*)+\\&& 2 Q_{1,2}\otimes
(iQ_{1,1}^*\wedge \wedge m_{13}^*)-3 Q_{1,2}\otimes
(iQ_{1,2}^*\wedge T_{1,1}^*)-2 iQ_{1,2}\otimes (Q_{1,1}^*\wedge
(T_{1,2}+T_{2,1})^*)+\\
&&2iQ_{1,2}\otimes (iQ_{1,1}^*\wedge i(T_{1,2}-T_{2,1})^*)+3 iQ_{1,2}\otimes (Q_{1,2}^*\wedge T_{1,1}^*)  \\[2mm]

c_6^2 & =&(A_{1,1}-A_{2,2})\otimes (Q_{1,1}^*\wedge
Q_{1,2}^*)-(A_{1,1}-A_{2,2})\otimes (iQ_{1,1}^*\wedge
iQ_{1,2}^*)+\\&& i(A_{1,1}-A_{2,2})\otimes (Q_{1,1}^*\wedge
iQ_{1,2}^*)+
i(A_{1,1}-A_{2,2})\otimes (iQ_{1,1}^*\wedge Q_{1,2}^*)+ \\
&&A_{1,2}\otimes (Q_{1,2}^*)^{\wedge 2}-A_{1,2}\otimes
(iQ_{1,2}^*)^{\wedge 2}
+ 2 i A_{1,2}\otimes (Q_{1,2}^*\wedge iQ_{1,2}^*)-A_{2,1}\otimes (Q_{1,1}^*)^{\wedge 2}+ \\
&& A_{2,1}\otimes (iQ_{1,1}^*)^{\wedge 2}- 2iA_{2,1}\otimes
(Q_{1,1}^*\wedge iQ_{1,1}^*)+ Q_{1,1}\otimes (Q_{1,1}^*\wedge
(T_{1,2}+T_{2,1})^*)-\\
&& Q_{1,1}\otimes (iQ_{1,1}^*\wedge
i(T_{1,2}-T_{2,1})^*)-Q_{1,1}\otimes (Q_{1,2}^*\wedge T_{1,1}^*)
-iQ_{1,1}\otimes (Q_{1,1}^*\wedge i(T_{1,2}-T_{2,1})^*)
-\\
&&iQ_{1,1}\otimes (iQ_{1,1}^*\wedge (T_{1,2}+T_{2,1})^*)+
iQ_{1,1}\otimes (iQ_{1,2}^*\wedge T_{1,1}^*)+Q_{1,2}\otimes (Q_{1,1}^*\wedge T_{2,2}^*) -\\
&&Q_{1,2}\otimes (Q_{1,2}^*\wedge (T_{1,2}+T_{2,1})^*)
-Q_{1,2}\otimes (iQ_{1,2}^*\wedge i(T_{1,2}-T_{2,1})^*)-iQ_{1,2}\otimes (iQ_{1,1}^*\wedge T_{2,2}^*)-\\
&&iQ_{1,2}\otimes (Q_{1,2}^*\wedge i(T_{1,2}-T_{2,1})^*)
+iQ_{1,2}\otimes (iQ_{1,2}^*\wedge (T_{1,2}+T_{2,1})^*)  \\[3mm]

c_7^2 & =& -(A_{1,1}-A_{2,2})\otimes (Q_{1,1}^*\wedge
iQ_{1,2}^*)-(A_{1,1}-A_{2,2})\otimes (iQ_{1,1}^*\wedge
Q_{1,2}^*)+\\&&i(A_{1,1}-A_{2,2})\otimes (Q_{1,1}^*\wedge
Q_{1,2}^*)- i(A_{1,1}-A_{2,2})\otimes (iQ_{1,1}^*\wedge
iQ_{1,2}^*)- 2A_{1,2}\otimes (Q_{1,2}^*\wedge iQ_{1,2}^*)+\\
&&i A_{1,2}\otimes (Q_{1,2}^*)^{\wedge 2}-i A_{1,2}\otimes
(iQ_{1,2}^*)^{\wedge 2}+2 A_{2,1}\otimes (Q_{1,1}^*\wedge
iQ_{1,1}^*)- iA_{2,1}\otimes (Q_{1,1}^*)^{\wedge
2}+\\&&iA_{2,1}\otimes (iQ_{1,1}^*)^{\wedge 2}-Q_{1,1}\otimes
(Q_{1,1}^*\wedge i(T_{1,2}-T_{2,1})^*)-Q_{1,1}\otimes
(iQ_{1,1}^*\wedge (T_{1,2}+T_{2,1})^*)+\\&& Q_{1,1}\otimes
(iQ_{1,2}^*\wedge T_{1,1}^*)- iQ_{1,1}\otimes (Q_{1,1}^*\wedge
(T_{1,2}+T_{2,1})^*)+ iQ_{1,1}\otimes (iQ_{1,1}^*\wedge
i(T_{1,2}-T_{2,1})^*)+ \\&& iQ_{1,1}\otimes (Q_{1,2}^*\wedge
T_{1,1}^*)-Q_{1,2}\otimes (iQ_{1,1}^*\wedge
T_{2,2}^*)-Q_{1,2}\otimes (Q_{1,2}^*\wedge
i(T_{1,2}-T_{2,1})^*)+\\&&Q_{1,2}\otimes (iQ_{1,2}^*\wedge
(T_{1,2}+T_{2,1})^*)-iQ_{1,2}\otimes (Q_{1,1}^*\wedge
T_{2,2}^*)+iQ_{1,2}\otimes (Q_{1,2}^*\wedge
(T_{1,2}+T_{2,1})^*)+\\
&&Q_{1,2}\otimes (iQ_{1,2}^*\wedge i(T_{1,2}-T_{2,1})^*)
\end{array}
\]
The two cocycles in degree $3$ become:

\[\tiny
\begin{array}{lcl}
c_8^3 & = & -A_{1,2}\otimes (Q_{1,1}^*\wedge
i(T_{1,2}-T_{2,1})^*)-A_{1,2}\otimes (iQ_{1,1}^*\wedge
(T_{1,2}+T_{2,1})^*)+A_{1,2}\otimes (iQ_{1,2}^*\wedge T_{1,1}^*)-\\
&& i A_{1,2}\otimes (Q_{1,1}^*\wedge (T_{1,2}+T_{2,1})^*)+i
A_{1,2}\otimes (iQ_{1,1}^*\wedge i(T_{1,2}-T_{2,1})^*)+i
A_{1,2}\otimes (Q_{1,2}^*\wedge T_{1,1}^*)+\\
&&Q_{1,2}\otimes (T_{1,1}^*\wedge i(T_{1,2}-T_{2,1})^*)+
iQ_{1,2}\otimes (T_{1,1}^*\wedge (T_{1,2}+T_{2,1})^*) \\[2mm]

c_9^3 & = &A_{1,2} \otimes (Q_{1,1}^*\wedge
(T_{1,2}+T_{2,1})^*)-A_{1,2}\otimes (iQ_{1,1}^*\wedge
i(T_{1,2}-T_{2,1})^*)-A_{1,2}\otimes (Q_{1,2}^*\wedge T_{1,1}^*)-\\
&&i A_{1,2}\otimes (Q_{1,1}^*\wedge i(T_{1,2}-T_{2,1})^*)- i
A_{1,2}\otimes (iQ_{1,1}^*\wedge
(T_{1,2}+T_{2,1})^*)+i A_{1,2}\otimes (iQ_{1,2}^*\wedge T_{1,1}^*)-\\
&&Q_{1,2}\otimes (T_{1,1}^*\wedge
(T_{1,2}+T_{2,1})^*)+iQ_{1,2}\otimes (T_{1,1}^*\wedge
i(T_{1,2}-T_{2,1})^*)
\end{array}
\]
\end{Theorem}

\parbegin{Comment}\label{comment2} All cocycles, except $c_5$, appear in pairs
that differ by multiplication by $i$, cf. Comment \ref{comment}, so
the irreducible module generated by any cocycle of the pair splits
into two irreducibles after complexification; the cocycle $c_5$
corresponds to the real $\fg_0$-module that remains irreducible
after complexification.
\end{Comment}

\sssec{The cases of $\cM_N$ for $N>1$ and various supermatrix
formats} Passing to $N>1$, it is convenient to complexify
$\fg_*=(\fg_{-1},\fg_0)_*$, as well as the $\fg_0$-module
$H^2(\fg_{-1};\fg_*)$, and consider in this  $\fg_0$-module only
highest weight vectors. We have $\fg_0^\Cee=\fs(\fgl_1(2)\oplus
\fgl(N)\oplus \fgl_2(2))$, where the index identifies a copy of
$\fgl(2;\Cee)$, whereas the operator $\fs(\cdot)$ singles out the
supertraceless part of the argument. Denote the elements of
$\fgl_2(2)$ by $C$ to distinguish them from the elements $A$ of
$\fgl_1(2)$ and use different names for $Q$ and $S:=``\overline
Q^t{}''$ which now are both complex and independent of each other,
compare with \eqref{M3}:
\begin{equation}\label{M31}
\begin{pmatrix}
A&V&U\\
Q&B&R\\
T&S&C
\end{pmatrix},\;\;\begin{array}{l}\text{where  $S, V\in\Mat_\Cee(2\times N)$,
$Q, R\in\Mat_\Cee(N\times 2)$,}\\
\text{$A, C\in\fgl(2;\Cee)$,
 $B\in\fgl(N)$, $\tr B=\tr(A+C)$.}\end{array}
\end{equation}
The elements $T$ are not hermitian now, as in \eqref{M1}, but
arbitrary elements of $\fgl(2;\Cee)$. Having found the highest
weight vectors of the complex representation of $\fg_0^\Cee$, we
have to check which of these irreducible representations are
complexifications of already complex representations of $\fg_0$
(these are to be found among those of multiplicity 2) and which are
of multiplicity 1 (complexifications of real representations), cf.
Comment \ref{comment}.

In \cite{GLs}, the calculations of $H^2(\fg_-; \fg_*)$ are performed
for $\fg_*=\fsl(4|N)$ realized in various supermatrix formats
($4|N$, $2|N|2$, and several other ones) for $N=1$, 2, 4 and 8, for
both conformal and reduced cases. Moreover, the calculations are
also performed for several types of parabolic subalgebras $\fp$,
{\bf smaller} (if $N>1$) than the one containing ``all above the
$Q$-and-$S$ diagonal''. These smaller subalgebras $\fp$ are chosen
so as to have the components of the non-holonomic curvature tensor
whose components independent on odd coordinates match those entering
the Einstein equations; for details, inessential in this paper, but
important, in our opinion, for understanding SUGRA, see \cite{GLs}.

\ssec{The integrability of almost complex structure of the
$1|2n$-dimensional over $\Cee$ supercurves with a distinguished
structure} When we study integrability of the almost complex
structure (preserved by the Lie (super)algebra $\faut(J)$) in
presence of some other structure (tensor or a distribution)
preserved (at the point) by the Lie (super)algebra $\fg_*$, we
should replace $\fg_*$ in formulas \eqref{curvtens} and
\eqref{Ncurvtens} by the (generalized) Cartan prolong of $(\fg_{-};
\fh_0)$, where $\fh_0:=\fg_0\cap \faut(J)$. Let $A+Bi\in\fg_0$ be
the decomposition into the real and imaginary part. Then the
supermatrices of $\fh_0:=\fg_0\cap \faut(J)$ are of the form
\eqref{ab-ba} with $A+Bi\in \fg_0$.

We consider superstrings with various additional structures. Let us
describe the Lie superalgebras over $\Cee$ that preserve,
infinitesimally, these additional structures. The general algebra
$\fvect(m|n)$ does not preserve anything, the divergence-free one
$\fsvect(m|n)$ preserves a volume element.

To both (centerless) $N$-extended Neveu-Schwarz and Ramond type
contact Lie superalgebras with Laurent polynomial coefficients only
one vectorial Lie superalgebra with polynomial coefficients
corresponds since locally the corresponding superstrings are
isomorphic. Similar is the case with the parametric family of
stringy superalgebras preserving a volume element. Recall the
description of vectorial Lie superalgebras with polynomial
coefficients over $\Cee$:

\[\fvect(m|n):=\fder\Cee[x,\theta]=\left\{\sum f_i\partial_{x_i}+\sum
g_j\partial_{\theta_j}\mid f_i, g_j\in\Cee[x,\theta]\right\},\]where
$x=(x_1,\dots, x_m)$ are even indeterminates,
$\theta=(\theta_1,\dots, \theta_n)$ are odd ones;
\[\fsvect(m|n):=\left\{\sum f_i\partial_{x_i}+\sum
g_j\partial_{\theta_j}\mid \sum \partial_{x_i}(f_i)+\sum
(-1)^{p(g_j)}
\partial_{\theta_j}(g_j)=0\right\};\]

$\fk(2n+1|m)$ preserves the distribution singled out by the odd form
\[\alpha_1=dt-\sum\limits_{i}(p_idq_i-q_idp_i)-\sum\limits_{j}(\xi_jd\eta_j+\eta_jd\xi_j)+
\begin{cases}0&\text{for $m$ even}\\
\theta d\theta&\text{for $m$ odd},\end{cases}\]where $t$, $p=(p_1,
\dots, p_n)$, $q=(q_1, \dots, q_n)$ are even indeterminates and
$\xi=(\xi_1, \dots, \xi_k)$, $\eta=(\eta_1, \dots, \eta_k)$ for
$2k=m$ for $m$ even (and $\theta$ for $m$ odd) are odd ones.
Set\index{$K_f$, contact vector field} \index{$H_f$, Hamiltonian
vector field}:
\begin{equation}
\label{2.3.1} K_f=(2-E)(f)\pder{t}-H_f + \pderf{f}{t} E \ \
\text{for any
$f\in\begin{cases}\Cee [t, p, q, \xi, \eta]&\text{for $m$ even}\\
\Cee [t, p, q, \xi, \eta, \theta]&\text{for $m$ odd},\end{cases}$}
\end{equation}
where $E=\sum\limits_i y_i \pder{y_{i}}$ (here the $y_{i}$ are all
the coordinates except $t$), and $H_f$ is the hamiltonian vector
field with Hamiltonian $f$:
\begin{equation}
\label{2.3.2'} H_f=\sum\limits_{i\leq n}\left (\pderf{f}{p_i}
\pder{q_i}-\pderf{f}{q_i} \pder{p_i}\right)
-(-1)^{p(f)}\sum\limits_{j\leq k}\left(\pderf{f}{\xi_j}
\pder{\eta_j}+ \pderf{f}{\eta_j} \pder{\xi_j}\right)+
\begin{cases}0&\text{for $m$ even}\\
\pderf{f}{\theta} \pder{\theta}&\text{for $m$ odd}.\end{cases}
\end{equation}

$\fm(n|n+1)$ preserves the distribution singled out by the even form
\[\alpha_0=d\tau+\sum\limits_{j}(\xi_jdq_j+q_jd\xi_j),\]where $q=(q_1, \dots, q_n)$
are even indeterminates and $\xi=(\xi_1, \dots, \xi_n)$, and $\tau$
are odd ones. For any $f\in\Cee [q, \xi, \tau]$, set:
\begin{equation}
\label{2.3.3} M_f=(2-E)(f)\pder{\tau}- Le_f -(-1)^{p(f)}
\pderf{f}{\tau} E, 
\end{equation}
where $E=\sum\limits_iy_i \pder{y_i}$ (here the $y_i$ are all the
coordinates except $\tau$), and
\begin{equation}
\label{2.3.4} Le_f=\sum\limits_{i\leq n}\left( \pderf{f}{q_i}\
\pder{\xi_i}+(-1)^{p(f)} \pderf{f}{\xi_i}\
\pder{q_i}\right ).
\end{equation}
\index{$M_f$, contact vector field} \index{$Le_f$, periplectic
vector field} Let $L_D$ be the Lie derivative along the vector field
$D$. Since
\begin{equation}
\label{2.3.5}
\renewcommand{\arraystretch}{1.4}
\begin{array}{l}
 L_{K_f}(\alpha_1)=2 \pderf{f}{t}\alpha_1=K_1(f)\alpha_1, \\
L_{M_f}(\alpha_0)=-(-1)^{p(f)}2 \pderf{
f}{\tau}\alpha_0=-(-1)^{p(f)}M_1(f)\alpha_0,
\end{array}
\end{equation}
it follows that $K_f\in \fk (2n+1|m)$ and $M_f\in \fm (n)$. It is
not difficult to show that \[\fk (2n+1|m)=\Span\left\{K_f\mid f\in\begin{cases}\Cee [t, p, q, \xi, \eta]&\text{for $m$ even}\\
\Cee [t, p, q, \xi, \eta, \theta]&\text{for $m$
odd}\end{cases}\right\},\]and $\fm (n)=\Span\{M_f\mid f\in \Cee [q,
\xi, \tau]\}$.

These vectorial Lie superalgebras $\fg_*=\oplus\fg_i$ are considered
with their \emph{standard} $\Zee$-grading in which the degree of
each indeterminate, except $t$ and $\tau$, is equal to 1, $\deg
t=\deg\tau=2$.

\sssbegin{Theorem}\label{string} The following are the simple
vectorial Lie superalgebras $\fg_*$ with polynomial coefficients in
their standard grading corresponding to the distinguished simple
stringy Lie superalgebras and the (highest weight with respect to
the $\fg_0$-action when appropriate) representatives of the basis
elements of $H^2(\fg_-; \fg_*)$:

$\bullet$ $\fvect(1|n)^\Ree\subset\fvect(2|2n;\Ree)$ for $n=1, 2$:
these cases are already considered in \eqref{n,0}.

$\bullet$ $\fsvect(1|2)^\Ree\subset\fvect(2|4;\Ree)$: The highest
weight cocycles representing the elements of $H^2(\fg_{-1}; \fg_*)$,
where $\fg_*=\fsvect(1|2)^\Ree$, are, in addition to those for
$\fvect(1|2)^\Ree$, as follows:
\begin{equation}\label{svect}
\footnotesize
\begin{array}{ll} c_3^1 &=\partial_1\otimes(
   \partial_4^*\wedge \partial_6^* - \partial_1^*\wedge
   \partial_3^*)+\partial_2 \otimes(
\partial_2^*\wedge\partial_3^*+\partial_5^*\wedge
   \partial_6^*)+\partial_3\otimes \left((\partial_3^*){}^{\wedge 2}+
   (\partial_6^*){}^{\wedge 2}\right)+\\
   &\partial_4\otimes
(\partial_1^*\wedge
   \partial_6^*-
   \partial_3^*\wedge \partial_4^*)+\partial_5\otimes(
\partial_3^*\wedge
   \partial_5^*-
   \partial_2^*\wedge \partial_6^*)\\

c_4^1&=\partial_1\otimes
\partial_3^*\wedge
   \partial_4^*-\partial_2\otimes
\partial_3^*\wedge \partial_5^*-\partial_3\otimes
   \partial_3^*\wedge \partial_6^*-\partial_4\otimes
   \partial_1^*\wedge \partial_3^*+\partial_5\otimes
\partial_2^*\wedge
   \partial_3^*+\partial_6\otimes
   (\partial_3^*){}^{\wedge 2}\\

c_5^1&=\partial_1\otimes
\partial_4^*\wedge \partial_6^*-\partial_2\otimes
\partial_5^*\wedge \partial_6^*-\partial_3\otimes
   (\partial_6^*){}^{\wedge
   2}-\partial_4\otimes
   \partial_1^*\wedge \partial_6^*+\partial_5\otimes
\partial_2^*\wedge
   \partial_6^*+\partial_6\otimes
   \partial_3^*\wedge
   \partial_6^*\\

c_6^1&=-\partial_1\otimes
\partial_1^*\wedge
   \partial_6^*+\partial_2\otimes
   \partial_2^*\wedge
   \partial_6^*+\partial_3\otimes
   \partial_3^*\wedge
   \partial_6^*-\partial_4\otimes
   \partial_4^*\wedge \partial_6^*+\partial_5\otimes
\partial_5^*\wedge
   \partial_6^*+\partial_6\otimes(
   \partial_6^*){}^{\wedge 2}
\end{array}
\end{equation}
\normalsize

$\bullet$ $\fk(1|1)^\Ree\subset\fvect(2|2;\Ree)$: We have
$\fg_-=\Span(K_1, iK_1, K_\theta, iK_\theta)$, $\fg_0=\Span(K_t,
iK_t)$ and ${H^2(\fg_{-}; \fg_*)=0}$.

$\bullet$ For $\fk(1|n)^\Ree\subset\fvect(2|2n;\Ree)$, where
$n\geq2$, all cocycles are of degree $0$:

$\bullet$ $\fk(1|2)^\Ree$: We have $\fg_-=\Span(K_1, iK_1, K_\xi,
iK_\xi, K_\eta, iK_\eta)$, $\fg_0=\Span(K_t, iK_t, K_{\xi\eta},
iK_{\xi\eta})$. cocycles of $H^2(\fg_{-}; \fg_*)$ are only the
following ones (hereafter, in cocycles, we write just $f$ instead of
$K_f$):
\[
\begin{array}{cclccl}
c_1^0&=&1\otimes \xi^*\wedge (i \xi)^*,&c_3^0&=&1\otimes (\xi^*)^{\wedge 2},\\

c_2^0&=&1\otimes \eta^*\wedge (i \eta)^*,&c_4^0&=& 1\otimes (\eta^*)^{\wedge 2 }\\

\end{array}
\]

$\bullet$ $\fk(1|3)^\Ree$: We have \[
\begin{array}{l}
\fg_-=\Span(K_1, iK_1, K_\xi, iK_\xi, K_\eta, iK_\eta, K_\theta,
iK_\theta),\\[2mm]
\fg_0=\Span(K_t, iK_t, K_{\xi\theta}, iK_{\xi\theta}, K_{\xi\eta},
iK_{\xi\eta}, K_{\eta\theta}, iK_{\eta\theta}).
\end{array}
\]
The highest weight cocycles of $H^2(\fg_{-}; \fg_*)$ with respect to
the $\fg_0$-action are
\[
c_1^0= 1\otimes (\eta^*)^{\wedge 2},\quad c_2^0= 1\otimes \eta^*
\wedge (i \eta)^*
\]

$\bullet$ For $\fk(1|2n)^\Ree$, where $n\geq2$, the highest weight
cocycles corresponding to the three irreducible components of the
$\fg_0$-module $H^2(\fg_{-}; \fg_*)$ are as follows:
\[
\tiny
\begin{array}{cclccl}
c_1^0 &=& 1\otimes (\xi_n^*)^{\wedge 2},&c_2^0&=& 1 \otimes \xi_n^* \wedge ( i \xi_n)^*,\\

c_3^0&=& 1 \otimes \xi_n^* \wedge \eta_1^*,
&c_4^0&=& 1 \otimes \xi_n^* \wedge ( i \eta_1)^*,\\

c_5^0&=& 1 \otimes (\eta_1^*)^{\wedge 2},& c_6^0&=& 1\otimes
\eta_1^* \wedge (i \eta_1)^*
\end{array}
\]

$\bullet$ For $\fk(1|2n+1)^\Ree$, where $n\geq2$, the highest weight
cocycles corresponding to the three irreducible components of the
$\fg_0$-module  $H^2(\fg_{-}; \fg_*)$ are as follows:
\[\tiny
\begin{array}{cclccl}
c_1^0&=& 1\otimes (\eta_1^*)^{\wedge 2},&c_2^0 &=& 1\otimes
(\eta_2^*)^{\wedge 2},\\

c_3^0&=& 1\otimes \eta_1^* \wedge (i \eta_1)^*,&c_4^0&=&  1\otimes
\eta_2^* \wedge (i \eta_2)^*,\\

c_5^0&=& 1\otimes \eta_1^* \wedge (\eta_2)^*,& c_6^0&=& 1\otimes
\eta_2^* \wedge (i \eta_2)^*.
\end{array}
\]

$\bullet$ $\fm(1)^\Ree\subset\fvect(2|4;\Ree)$: We have
\[\begin{array}{l}
\fg_-=\Span(M_1, iM_1, M_q, iM_q, M_\theta,
iM_\theta),\;\; 
\fg_0=\Span(M_\tau, iM_\tau, M_{q\theta}, iM_{q\theta}, M_{q^2},
iM_{q^2}).
\end{array}
\]
The cocycles representing the elements of $H^2(\fg_{-1}; \fg_*)$ are
as follows:
\[\tiny
\begin{array}{cclccl}
c_1^0 &=&  1\otimes   q_1 {}^*\wedge   (i \, q _1) {}^*,& c_3^1&=&
1\otimes   \theta _1{}^* \wedge
 1^*  + q_1\otimes    i \, q_1 {}^* \wedge   (i \, \theta
   _1) {}^*,\\
c_2^0&=& 1\otimes  (q_1{}^*) {}^{\wedge 2},& c_4^1&=&  1\otimes
\theta _1{}^* \wedge  (i \, 1)^*
    -q_1\otimes   q_1{}^* \wedge
(i \, \theta
   _1) {}^*;
\end{array}
\]
the operators $M_{q^2}$ and $iM_{q^2}$ trivially act on the
cohomology classes represented by the above cocycles.
\end{Theorem}

\parbegin{Remark} Although $\fvect(1|1)\cong\fk(1|2)\cong\fm(1)$ as
abstract Lie superalgebras, see \cite{GLS}, they are non-isomorphic
as filtered or $\Zee$-graded ones, and preserve completely different
structures. So no wonder that the obstructions to integrability to
the almost complex (resp. almost real-complex) structure (in the
first case, resp. the other two cases) look  (and are) completely
different.\end{Remark}

\section{On K\"ahler and hyper-K\"ahler supermanifolds}

\ssec{Definition on manifolds} Let a real manifold $M$ possess an
almost complex structure $J$ and a non-degenerate symmetric bilinear
form $h$ such that
\begin{equation}\label{Herm}
h(X, Y)=h(JX, JY) \text{~~for any vector fields~~} X, Y\in\fvect(M)
\end{equation} (such $h$ is said to be \emph{pseudo-Hermitian}). The
manifold $M$ is said to be \emph{K\"ahler} if $J$ is covariantly
constant with regard to the Levi-Civita connection $\nabla$
corresponding to the metric $h$, i.e.,
\begin{equation}\label{nab}
\nabla J=0.
\end{equation}
 Each K\"ahler manifold is symplectic in a
natural way with the non-degenerate 2-form defined by
\begin{equation}\label{symp}
\omega(X, Y)=h(JX, Y) \text{ for any $X, Y\in\fvect(M)$},
\end{equation}
whereas requirement
\begin{equation}\label{domega}
d\omega=0
\end{equation}
is one of definitions of K\"ahler manifolds instead of \eqref{nab}.

Any two of the constituents of the triple $(\omega, h, J)$ determine
the third one by means of eq. \eqref{symp} and, since on
supermanifolds these two entities can be even or odd, the notion of
K\"ahler manifold has (at least) four types of superizations.

\sssbegin{Remark} For the formula \eqref{symp} to define any third
ingredient of the triple $(\omega, h, J)$ given the other two, we
only need non-degeneracy of $\omega$ and $h$. The flatness of the
$G$-structures associated with  $J$ in the conventional definition
of the K\"ahler manifold (i.e., requirements  that $J$ is complex,
not almost complex) does not seem to be justified: We do not require
flatness of the metric, so why discriminate $\omega$ and $J$?
Besides, why should the metric be sign-definite?

M.Verbitsky informed us that, indeed, the sign-definiteness of $h$
in the traditional definition of the K\"ahler manifold is
unnecessary (published classification results are only known,
however, for sign-definite forms $h$), whereas the symplectic
structure is needed because
\begin{equation}\label{domnab} d\omega=0\Longrightarrow \nabla J=0.
\end{equation}
The condition $d\omega=0$ is, actually, a system of two equations
(vanishing of both components of
$d\omega\in\Omega^3=P^3\oplus\omega\wedge\Omega^1$, where $P^3$ is
the space of primitive (aka harmonic) forms, the ones ``not
divisible'' by $\omega$), and both are needed to ensure \eqref{nab}.

Moreover, the condition $d\omega=0$ is often taken for the
definition of the K\"ahler manifold.
\end{Remark}

In view of this Remark, we suggest the following definition of
K\"ahler supermanifold suitable also in the non-super setting.

\ssec{Definitions on supermanifolds} A~non-degenerate supersymmetric
bilinear form $h$ on the superspace $V$ will be called
\emph{pseudo-hermitian metric} relative the operator $J\in\End(V)$
such that $J^2=\pm\id$ if
\begin{equation}\label{sHerm}
h(X, Y)=(-1)^{p(X)p(J)}h(JX, JY) \text{~~for any vectors~~} X, Y\in
V.
\end{equation}
Let $\cM$ be a~real supermanifold with a complex structure (or a
$J$-symmetric or $\Pi$-symmetric structure over any ground field),
$h$ a~non-degenerate pseudo-hermitian metric relative to $J$ or
$\Pi$, and $\omega$ a~non-degenerate differential 2-form (for
details of definition of these notions, see \cite{L1}). The
supermanifold $\cM$ is said to be {\it K\"ahler} (an {\it almost}
one if $J$-symmetric or $\Pi$-symmetric or the complex structure is
an \lq\lq almost" one) \index{supermanifold, almost K\"ahler} if
\begin{equation}
\label{eq1}
\renewcommand{\arraystretch}{1.4}
\begin{array}{l}
\omega(X, Y)=h(JX, Y) \; \text{ for any $X,
Y\in\fvect(\cM)$ provided $p(h)+p(J)=p(\omega)$ and}\\
\nabla J=0\; \text{ for the Levi-Civita connection $\nabla$
corresponding to the metric $h$, see \cite{Po}.} \end{array}
\end{equation}
This definition implies the following restrictions on the possible
superdimensions of $\cM$:
\begin{equation}
\label{eq2}
\renewcommand{\arraystretch}{1.4}
\begin{array}{|c|c|c|}
\hline &p(J)=\ev&p(J)\text{~~(or $p(\Pi)$)}=\od\cr \hline
p(h)=\ev&2n|2m&2n|2n\cr \hline p(h)=\od&2n|2n&n|n\cr \hline
\end{array}
\end{equation}
The integrability conditions are not, however, automatically
satisfied and must be verified. (Recall that the odd non-degenerate
and closed form $\omega$ obtained in the off-diagonal cases of table
\eqref{eq2} is called, as A.Weil suggested, \emph{periplectic}, it
is the analog of the symplectic form that leads to the
``antibracket'' and the corresponding mechanics, discovered in
\cite{Lm} and rediscovered by I.Batalin and G.Vilkovissky who
applied it to physics \cite{BV}; for a review, see \cite{GPS}.)

The Lie superalgebra preserving an odd non-degenerate anti-symmetric
bilinear form is called \emph{periplectic} and denoted $\fpe^a(n)$;
an odd non-degenerate symmetric bilinear form is preserved by a Lie
superalgebra $\fpe(n)$, isomorphic to $\fpe^a(n)$, but having a
different matrix realization, see \cite{L1}.

Given three (almost) complex structures $J_i$ satisfying the
relations of quoternionic units, and one metric $h$ pseudo-hermitian
relative each $J_i$, together with three symplectic (or periplectic)
forms $\omega_i$ tied together by three relations of the form
\eqref{eq1}, we arrive at the notion of an (almost)
\emph{hyper-K\"ahler supermanifold}.

Most of the notions of this section had never been distinguished
before. It would be interesting to generalize with their help the
ideas exposed, e.g., in lectures \cite{Ku} and later works.

\section{Appendix: Cartan prolongations and its generalizations (\cite{Shch})}

\subsection{Cartan prolongations}\label{s15.2.5.0}
Let $\fg$ be a Lie algebra, $V$ a $\fg$-module, and $S^i$ the
operator of the $i$th symmetric power. Set $\fg_{-1} = V$, $\fg_0 =
\fg$, and for $k
> 0$, define the $k$th \emph{Cartan prolongation} of the pair
$(\fg_{-1}, \fg_0)$ by setting
\begin{equation*}
\renewcommand{\arraystretch}{1.4}
\begin{split}
    \fg_k= \{X\in \Hom(\fg_{-1}, \fg_{k-1})\mid & X(v_0)(v_1, v_2,\dots v_k) =
  X(v_1)(v_0, v_2, \dots, v_k)\\
 &\text{ for any }\; v_0, v_1, \dots v_k\in \fg_{-1}\}.
\end{split}
\end{equation*}
Let
\begin{equation}
 \label{clsfeq120}
\begin{gathered}
 i\colon S^{k+1}(\fg_{-1})^*\otimes \fg_{-1}\tto
 S^{k}(\fg_{-1})^*\otimes \fg_{-1}^*\otimes\fg_{-1},\\
 j\colon S^{k}(\fg_{-1})^*\otimes \fg_{0}\tto
 S^{k}(\fg_{-1})^*\otimes \fg_{-1}^*\otimes\fg_{-1}
 \end{gathered}
\end{equation}
be the natural embeddings. Then
{\MathSkip{.5}$\fg_k=i(S^{k+1}(\fg_{-1})^*\otimes \fg_{-1})\cap
j(S^{k}(\fg_{-1})^*\otimes \fg_{0})$}.

\emph{The complete Cartan prolong} of the pair $(V, \fg)$ is the
space $(\fg_{-1},\fg_{0})_* = \boplus_{k\geq -1} \fg_k$. This space
is naturally endowed with a Lie algebra structure which is rather
bothersome to define in abstract terms. If, however, the
$\fg_0$-module $\fg_{-1}$ is \emph{faithful}, there is an embedding
\[
\begin{gathered}
 (\fg_{-1}, \fg_{0})_*\subset \fvect (n) =
 \fder \Cee[x_1,\ldots , x_n],\ \text{ where }\ n = \dim \fg_{-1}\ \text{ and }\\
 \fg_i = \{D\in \fvect(n)\mid \deg D=i, [D, X]\in\fg_{i-1}\ \text{ for any }
 X\in\fg_{-1}\},
\end{gathered}
\]
and the Lie algebra structure on $\fvect (n)$ induces same on
$(\fg_{-1}, \fg_{0})_*$, the latter structure coincides with the one
we were lazy to define in abstract terms.

Of four series of simple vectorial Lie algebras with polynomial
coefficients, three are the complete Cartan prolongs:
\begin{equation*}
 \fvect(n)=(\id, \fgl(n))_*, \quad \fsvect(n)=(\id,
\fsl(n))_*,\quad \fh(2n)=(\id, \fsp(n))_*,
\end{equation*}
The fourth series ---  $\fk(2n+1)$ --- is the result of a bit more
general construction to be described in subsec.\,\ref{s15.2.5.1}.

\subsubsection{Vectorial Lie superalgebras as Cartan prolongs}\label{s15.2.5.2}
Superization of constructions of subsec.\,\ref{s15.2.5.0} is direct
one: via the Sign Rule. We obtain in this way the following infinite
dimensional Lie superalgebras:
\begin{equation*}
\begin{split}
 &\fvect(m|n)=(\id, \fgl(m|n))_*, \quad \fsvect(m|n)=(\id, \fsl(m|n))_*; \\
 &\fh(2m|n)=(\id, \fosp^{a}(m|2n))_*, \quad \fle(n)=(\id, \fpe^{a}(n))_*,\\
\end{split}
\end{equation*}
where $\fosp^{a}(m|2n)$ is assumed preserving a non-degenerate even
anti-symmetric bilinear form, while an isomorphic to it Lie
superalgebra $\fosp(m|2n)$ is assumed preserving a non-degenerate
even symmetric bilinear form.

\begin{Remarks}
1) The Cartan prolong $(\id,
\fosp(m|2n))_*=(\Pi(\id),\fosp^a(m|2n))_*$ is of finite dimension.

2) Superization of the contact Lie algebras leads to the {\bf two}
series, $\fk$ and $\fm$, described below.
\end{Remarks}

\subsection{Generalizations of Cartan prolongation}\label{s15.2.5.1}
Consider a nilpotent $\Zee$-graded Lie algebra
$\fg_-=\boplus_{-d\leq i\leq -1}\fg_i$, and a subalgebra
$\fg_0\subset \fder_0\fg$ of the Lie algebra of its $\Zee$-grading
preserving derivations. Let
\begin{equation*}
\begin{gathered}
 i\colon S^{k+1}(\fg_{-})^*\otimes \fg_{-}\tto
 S^{k}(\fg_{-})^*\otimes \fg_{-}^*\otimes\fg_{-}, \\
 j\colon S^{k}(\fg_{-})^*\otimes \fg_{0}\tto
 S^{k}(\fg_{-})^*\otimes \fg_{-}^*\otimes\fg_{-}
\end{gathered}
\end{equation*}
be narural embeddings analogous to \eqref{clsfeq120}. For $k>0$,
define the $k$th prolongation of the pair $(\fg_-, \fg_0)$ by
setting
\begin{equation*}
\fg_k = \left (j(S^{\bcdot}(\fg_-)^*\otimes \fg_0)\cap
i(S^{\bcdot}(\fg_-)^*\otimes \fg_-)\right )_k,
\end{equation*}
where the index $k$ in the right hand side singles out a component
of degree $k$.

Set $(\fg_-, \fg_0)_*=\mathop{\oplus}\limits_{i\geq -d} \fg_i$. If
$\fg_0$-modules $\fg_i$ are irreducible for $i<0$, then, as is easy
to check, $(\fg_-,\fg_0)_*$ is a Lie subalgebra in
$\fvect(\dim\fg_-)$.

Superization of the construction is obvious.

What is the Lie algebra of contact vector fields in these terms? Let
$\fhei(2n)$ be the Heisenberg Lie algebra: its space is $W\oplus
{\Cee}\cdot z$, where $W$ is a $2n$-dimensional space endowed with a
non-degenerate anti-symmetric bilinear form $B$, and the bracket is
given by the following relations:
\begin{equation}
\label{clsfeq126} \text{$z$ lies in the center and $[v, w]=B(v,
w)\cdot z$ for any $v, w\in W$}.
\end{equation}

Clearly, $ \fk(2n+1)\cong (\fhei(2n), \fc\fsp(2n))_*$, where
$\fc\fg$ denotes the trivial central extension of the Lie algebra
$\fg$ by means of a 1-dimensional center.

$\bullet$ Define a structure of a Lie superalgebra $\fhei(2n|m)$ on
the direct sum of a $(2n|m)$-dimensional superspace $W$ endowed with
a non-degenerate {\bf even} anti-symmetric bilinear form $B$ and a
$(1, 0)$-dimensional space spanned by a vector $z$ by the expression
\eqref{clsfeq126}. Obviously,
\begin{equation*}
 \fk(2n+1|m)=(\fhei(2n|m), \fc\fosp^{a}(m|2n))_*.
\end{equation*}

$\bullet$ Pericontact (i.e., \lq\lq odd" contact) analog of the
serias $\fk$ is associated with the following \lq\lq odd" analog of
the Lie superalgebra $\fhei(2n|m)$. Let $\fab(n)$ be the
\emph{antibracket} superalgebra: its space is $W\oplus \Cee\cdot z$,
where $W$ is an $n|n$-dimensional superspace endowed with a
non-degenerate {\bf odd} anti-symmetric bilinear form  $B$, and the
bracket in $\fab(n)$ is given by the following relations:
\begin{equation*}
 \text{$z$ is odd and lies in the center, and $[v,
 w]=B(v, w)\cdot z$ for any $v, w\in W$.}
\end{equation*}

Clearly,
\begin{equation*}
 \fm(n)=(\fab(n), \fc\fpe^{a}(n))_*.
\end{equation*}


\end{document}